\documentclass{notices}
\usepackage{amsfonts,amssymb,amsmath,amscd,graphicx}

\def\ve{\varepsilon}
\def\half{\textstyle{1\over 2}}
\def\Li{\hbox{Li}}
\def\Re{\hbox{Re}\kern 1pt }

      \def\zoo{1}
    \def\apple{2}
     \def\bomb{3}
  \def\chenhou{4}
      \def\cmi{5}
     \def\cook{6}
       \def\ct{7}
  \def\cooktsp{8}
     \def\feff{9}
     \def\fort{10}
  \def\gourdon{11}
     \def\mont{12}
       \def\pt{13}
   \def\sarnak{14}
 \def\yogiisms{15}
  \def\pinball{16}
     \def\wang{17}
      \def\npc{18}
     \def\will{19}

\title{Reflections on the Millennium Problems}

\author{
Lloyd N. Trefethen
\affil{Professor of Applied Mathematics in Residence,
School of Engineering and Applied Sciences, Harvard University}
}

\begin{document}
\date{}

\maketitle

Mathematics is the only discipline that can pose a problem at one point in time
and unambiguously solve it a century later.
It is an astonishing tribute
to the rigor and abstract nature of the field that we can do this.
There are famous examples such as Fermat's Last Theorem (conjectured in
1637, proved in 1994) and the Poincar\'e conjecture (1904, 2002)---and
innumerable other instances on at least decadal time scales.
Theorems including the infinitude of the primes and the irrationality
of $\sqrt 2$ go back to the ancient Greeks, and the proofs remain
valid today.  The twin prime conjecture dates to at least 1849
and remains open, but such is our confidence in the longevity of
mathematics that a recent Fields Medal recognized work that may be
a step toward its eventual resolution.  Other disciplines, whether
physics or philosophy or histology or history, are not like this.

\begin{figure}
\begin{center}
\includegraphics[clip,scale=.27]{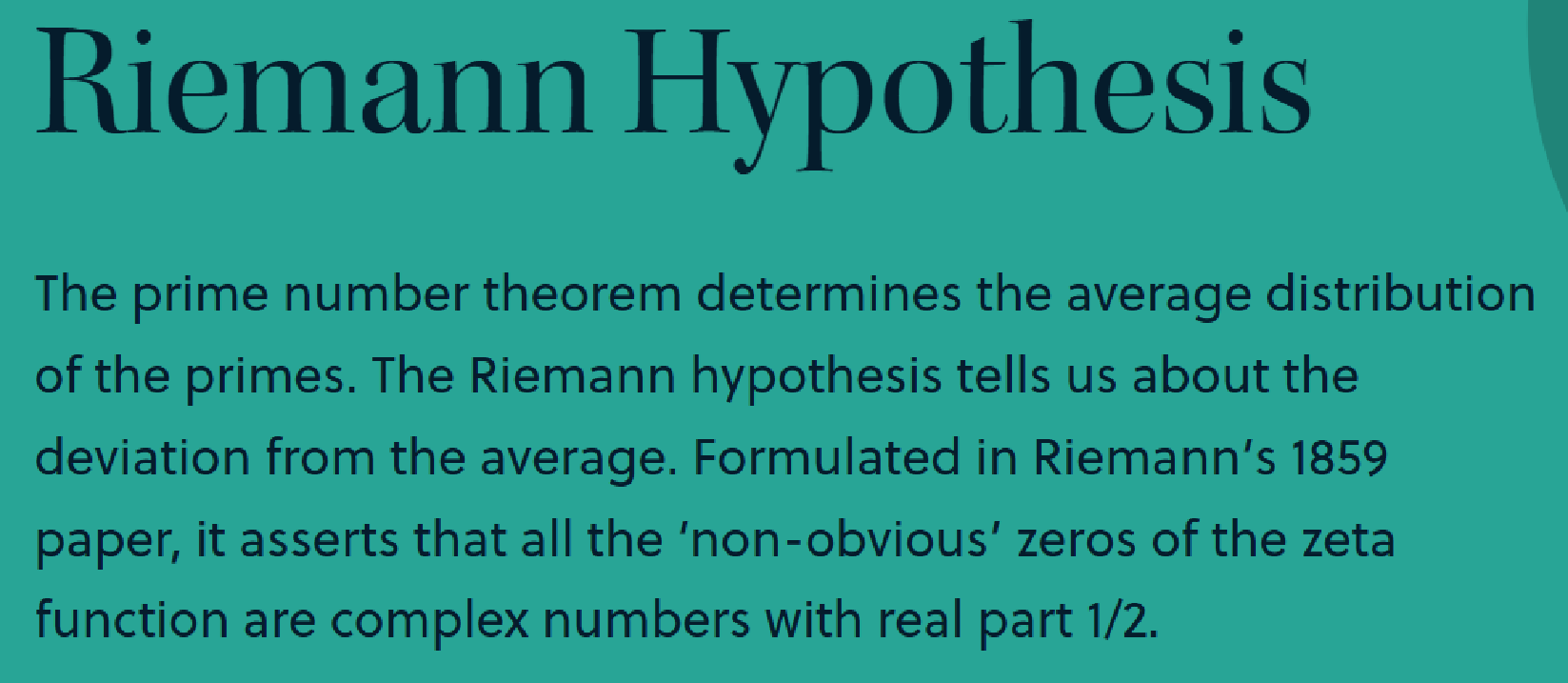}
\vspace*{6pt}

\includegraphics[clip,scale=.27]{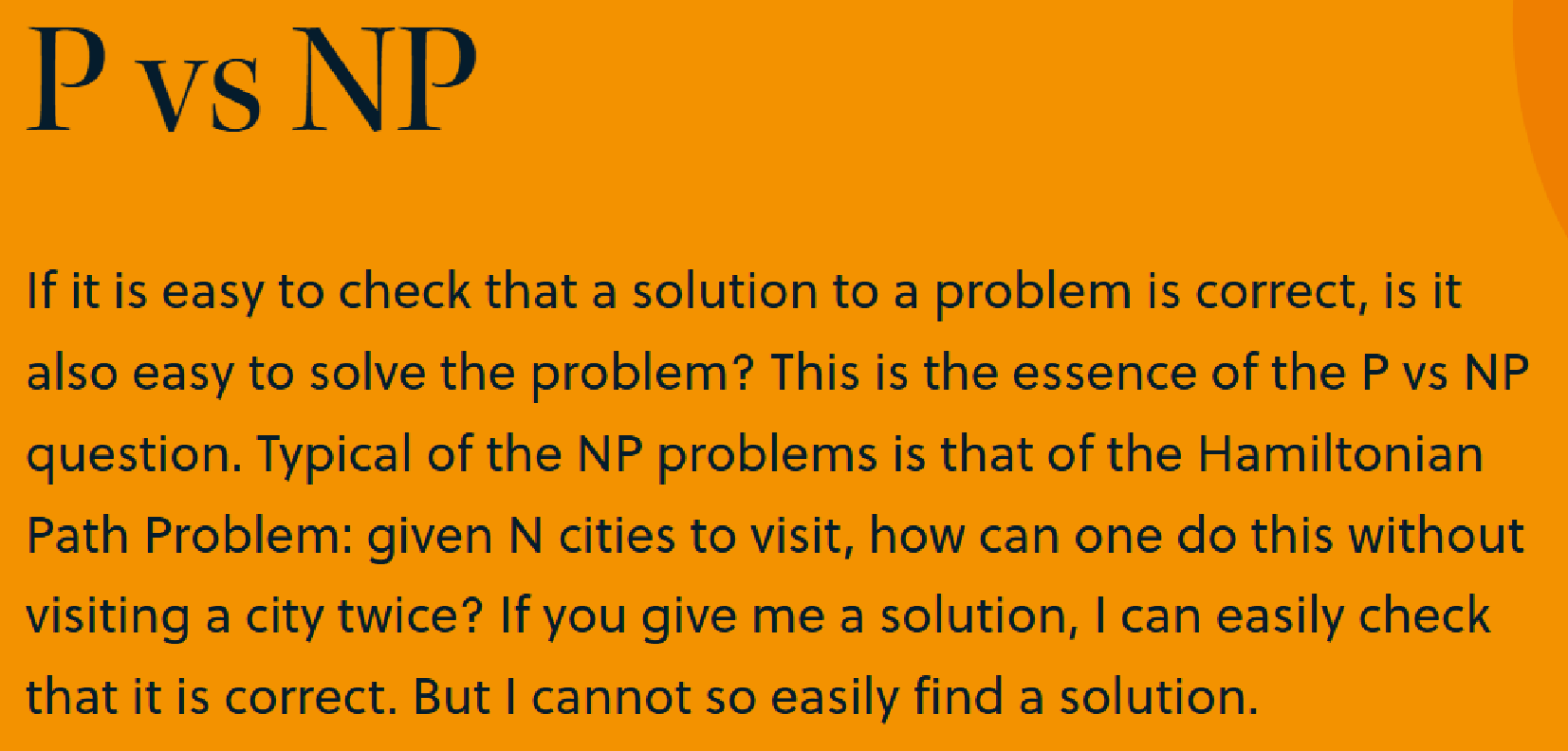}
\vspace*{7pt}

\includegraphics[clip,scale=.277]{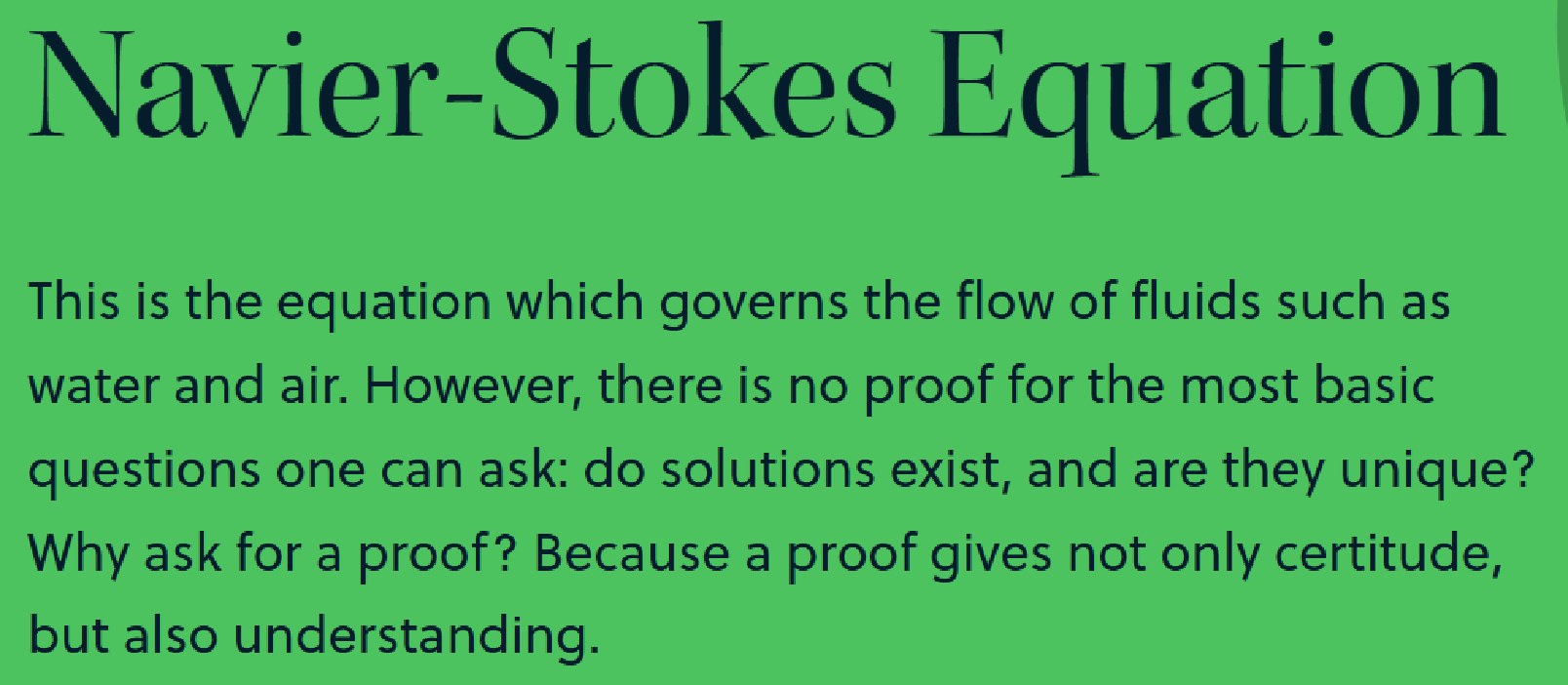}
\end{center}
\vspace*{-8pt}
\caption{\label{three}Top-level summaries of the three Prize problems
at the CMI web site [\cmi].}
\end{figure}

Seven famous examples are the \$1-million Millennium Prize Problems announced
by the Clay Mathematics Institute (CMI) in 2000 [\cmi].  Each can be
stated loosely in a few sentences (Figure~\ref{three}) and was given
a precise presentation by an expert in the area [\bomb,\cook,\feff].
The Poincar\'e conjecture has been settled, but the other 
Millennium problems remain open.  In this note, I reflect on
three of these: the Riemann Hypothesis, P vs.\ NP,
and the solvability of the Navier-Stokes equations.

The details differ from problem to problem, which makes it all the
more remarkable that a common theme emerges.  Each of these three
problems was born as a challenge whose importance, at least in part,
was due to its potential quantitative or practical implications.
However, as time has passed and the problem has eluded solution, each
has lost a good deal of its original leverage---for three curiously
different reasons.  Quietly over the years, mostly without comment
and sometimes without much recognition by experts or non-experts,
all three have consolidated into primarily academic challenges, whose resolution
will still have great impact---but not so much in the directions initially
envisioned among mathematicians and presented to the world.

\section{Riemann Hypothesis}
The Prime Number Theorem asserts that
\begin{equation}
\pi(x) \sim \Li(x) = \int_2^x {dt\over \log t}\kern 2pt,
\label{pnt0}
\end{equation}
where $\pi(x)$ denotes the number of primes less than or equal to $x$.
The logarithmic integral is asymptotic to
$x/\log x$, and thus for any $\ve>0$, $\pi(x)$ eventually exceeds $x^{1-\ve}$.
The original question behind the Riemann Hypothesis (RH)
is, how big is the error in (1)?  It is known that
if the Riemann zeta function $\zeta(s)$ has a zero with real part $\beta> \half$,
then
\begin{equation}
|\pi(x) - \Li(x)| > {C\kern .5pt x^\beta\over \log x}
\label{pnt}
\end{equation}
for some $C>0$ for an unbounded set of values of $x$.
The hypothesis asserts that all the complex zeros
of $\zeta(s)$ lie on the critical line $\Re s = \half$. 
It is known that RH is true if and only if
\begin{equation}
\pi(x) - \Li(x) = O\bigl(x^{1/2} \log x\bigr),
\end{equation}
whereas if it is false, then there is
some $\beta\in (\half,1) $ such that
(\ref{pnt}) holds for an unbounded set of $x$ values.
(Note that any zeros off the line would have to
come in pairs with real parts $\beta$ and $1-\beta$.)

It was Riemann in 1859 who brilliantly investigated the behavior of
$\zeta(s)$ in the complex $s$-plane,
calculated that the first few zeros
were on the critical line, and hypothesized that this might be true
for all of them.
His paper is called ``On the number of primes less than a given size.''
(A key predecessor was Euler.)
Over the nearly 170 years since, the interest of number theorists in
this problem has been enormous.
Here are some words of Bombieri in his statement of the problem
in 2000 [\bomb]:
\begin{quote}
\noindent The failure of the Riemann hypothesis would create havoc in
the distribution of prime numbers.
\end{quote}

\begin{figure}
\begin{center}
\includegraphics[clip,scale=.5]{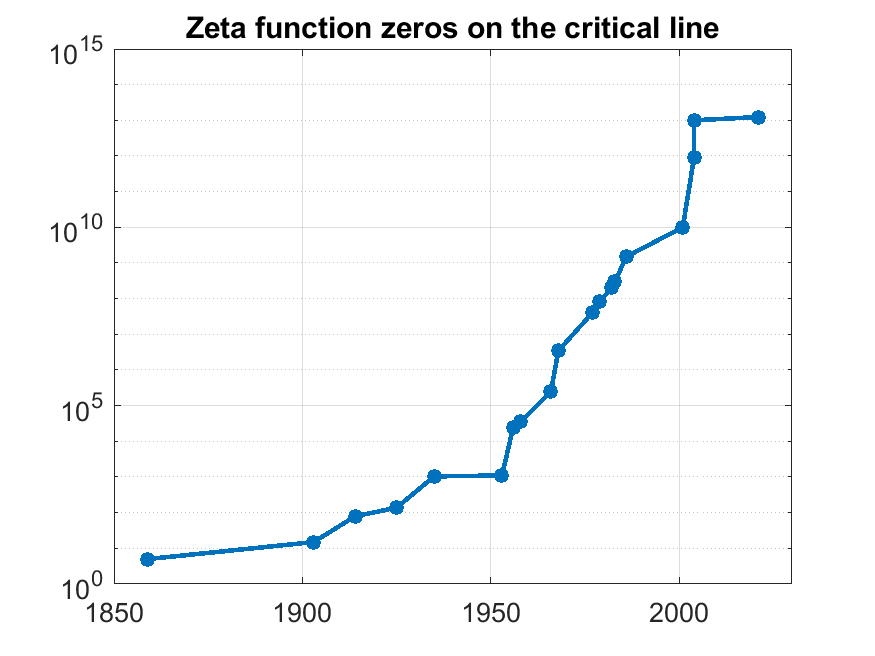}
\vspace*{-18pt}
\end{center}
\caption{\label{zeros}Number of zeros of $\zeta(s)$ known to lie on the critical line.
The increase from $10$ to $10^{13}$ reduces 
the quantitative implication of a potential failure of the Riemann Hypothesis by a
factor on the order of a trillion.
Data from [\gourdon] and [\pt].}
\end{figure}

Yet it is interesting to ask: What would actually be the quantitative
implications of failure of RH in terms of (1)--(3)?  Here is where
170 years of history come into play.  Riemann calculated three
of the zeros, or possibly four or five.  (His calculations
were not rediscovered until the 1930\kern .4pt s.)  In the years
since, using methods that owe a good deal to Turing,
an astonishing {\em twelve trillion\/} zeros have been
computed and found to lie on the line (Figure~\ref{zeros}) [\pt].
(The exact number is $12{,}363{,}153{,}437{,}138$.)

The consequence of these calculations goes
to the constant $C$ of (\ref{pnt}).
Suppose RH is proved to be false tomorrow, so that some zero beyond the
first $10^{13}$ has real part $\beta > \half$.
Roughly speaking (not stated rigorously), this will imply that (\ref{pnt}) becomes
\begin{equation}
|\pi(x) - \Li(x)| > {10^{-13}\kern .5pt x^\beta\over \log x}.
\label{pnt2}
\end{equation}
(For precise estimates, see [\mont] and [\pt].)  Such a lower
bound would be extraordinarily small in comparison to $\pi(x)$.
I believe it would be so small that the effect of failure of RH on
$\pi(x)$ would be undetectable with current computing technology.

So, is the Riemann Hypothesis $99.99999999999\%$ proved?  Of course
not.  Its proof or disproof will still be an earthquake, when
and if it comes.  But this is not because of what it would 
teach us about the distribution of primes.  More 
exciting are the rich new directions in which the subject has developed
since Riemann's time.  Hadamard and de La Vall\'ee Poussin finally proved
the Prime Number Theorem in 1896
by showing that $\zeta(s)$ has no zeros with $\Re s = 1$ (a
19\kern .2pt th-century-long open problem eminently worthy of a prize!), and
RH has been generalized to a broad class of
{\em $L$-functions\/}, all of which may have their zeros on the line
$\Re s = \half$.  One reason why experts would be
so pleased to see a proof of RH is that it might
introduce methods leading to a proof of this Generalized Riemann Hypothesis (GRH), from
which many other results would follow.  Peter
Sarnak put it like this in an April 2026 lecture at Harvard~[\sarnak]:
\begin{quote}
There aren't any
knockout consequences of the Riemann Hypothesis
that I know of$\kern 1pt\dots.$ It's an interesting problem when you add
``G'' to ``RH''.
\end{quote}

\section{P vs.\ NP}
In my perhaps simplified understanding of the history of computing,
something changed in the 1960\kern .4pt s, around the time of the
introduction of the Fast Fourier Transform (FFT).  The FFT speeds up
$n$-point transforms from $O(n^2)$ to $O(n\log n)$ operations, and
when the paper by Cooley and Tukey appeared in 1965 [\ct], scientists and
engineers saw that differences in asymptotic operation counts would
be decisive as machines got faster.  (Cooley and Tukey reported a
computing time of 8 seconds for an FFT of length 8192; it takes
20 microseconds on my laptop.)  Moore's Law also dates to 1965,
and the future seemed clear.  Problems would get bigger on faster
machines, and the gaps between $O(n\log n)$, $O(n^2\kern 1pt )$,
and $O(n^3\kern 1pt )$ would grow.  Exponential time algorithms,
requiring $O(C^n)$ operations for some $C>1$, would be left far
behind.

It wasn't long after this the most famous problem of theoretical
computer science was formulated.  
Around 1971 Cook, Levin, and Karp discovered that a large class
of {\em NP-complete\/} problems could only be
solved as far as anyone knew by exponential algorithms, yet were
all equivalent in the sense that any one could be reduced to any
other with only a polynomial amount of work.  (Key predecessors were
G\"odel and Edmonds.)  Thus there were two
great classes of problems.  Problems in P, such as
inverting a matrix or sorting, could be solved in polynomial
time, whereas NP-complete problems, such as graph coloring or the Traveling
Salesman Problem, were all equivalent and probably hard.
Possibly a polynomial
algorithm for an NP-complete problem would be discovered, and then
P and NP would be the same.  More likely they were different, and the challenge
was to prove this.  What computational question
could be more consequential?  If P and NP were different,
the gulf between the two classes was destined to widen inexorably.
A big part of the future of computing seemed to be tied to the
problem of P and NP.
\begin{figure}
\begin{center}
\includegraphics[clip,scale=.5]{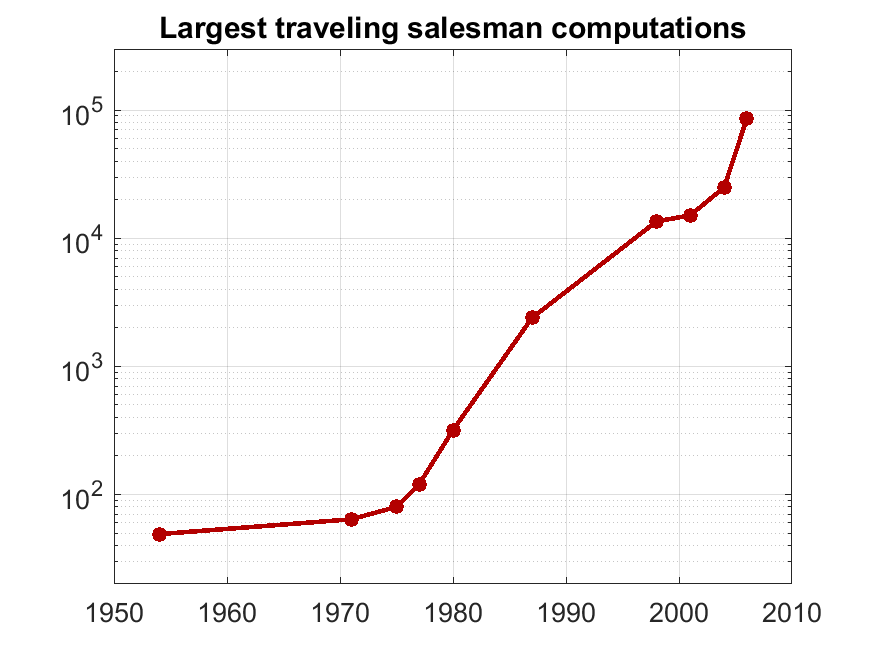}
\vspace*{-18pt}
\end{center}
\caption{\label{tspfig}Traveling Salesman Problem computations of record size.
Despite the exponential worst-case complexity of known algorithms,
the scale has increased 1000-fold since the P vs.\ NP problem was posed in
the early 1970\kern .5pt s.  Data from [\apple] and [\cooktsp].}
\end{figure}

This is not how history has unfolded.\footnote{I have discussed P
and NP previously in [\yogiisms].} Here in 2026, with computers a
million times faster, polynomial and exponential algorithms should
have long since lost sight of each other.  Instead, we see a wide
spectrum of behaviors.  Problems with polynomial algorithms are
indeed mostly easy enough to solve.  Problems whose only known
algorithms are exponential, on the other hand, are all over the
place.  Some seem to be truly hard, but others are surprisingly
tractable (Figure~\ref{tspfig}).

Two phenomena explain this effect.  The first is that some algorithms
that take exponential time in the worst case often run much faster, and
algorithms have been tuned to make the most of such gaps.  For an example
one need look no further than the simplex algorithm for linear
programming, which is provably exponential in the worst case.
The simplex algorithm is so fast and reliable that it remains
a workhorse of large-scale optimization, even though linear
programming belongs to P and polynomial algorithms are known.
There are also many instances of far-from-worst-case behavior among
the NP-complete problems.  A conspicuous example is the 3-SAT problem
of testing satisfiability of large Boolean expressions in
3-conjunctive normal form, one of the original NP-complete problems
established by Cook and Levin.  In 1971, who would have predicted
that 3-SAT solvers would become a trusted tool for industrial
applications?  I remember my surprise when Don Knuth, the creator
of the field of analysis of algorithms, gave the SIAM von Neumann
Prize lecture in 2016 on ``Satisfiability and combinatorics.''
The lecture was all about 3-SAT applications, yet Knuth mentioned
only in passing that the problem was NP-complete.

The second development that has blurred the boundary between P
and NP is that there are many problems where exponential behavior
is eliminated even in the worst case, typically by exploiting
numerical methods of continuous optimization, provided one is
willing to settle for somewhat less than optimality.  For example,
this is the case with max cut and various other graph problems,
including the Euclidean Traveling Salesman Problem.  In the case
of max cut, a polynomial algorithm can get within 88\%
of the optimum [\will].

The possibility is still alive that a breakthrough algorithm will
be discovered that solves an NP-complete problem in fast polynomial
time, and that would be a bombshell.  But
hopes have faded that this is likely to happen.

So, P vs.\ NP is not the problem it seemed 50 years ago, but
paradoxically, its visibility as a pillar of theoretical computer
science has only grown.  Thousands of NP-complete problems have
been identified [\npc], and according to the Complexity Zoo as
of August 2026, at least 551 different complexity classes have
been investigated, from A${}_0$PP to ZQP [\zoo].\ \ The sense
of P vs.\ NP as a theoretical organizing principle is stronger than ever,
and innumerable results have been developed in this framework.
There is also fascination with the idea that if P were equal to NP,
then much of mathematics might in theory crumble away, since solving
hard problems---including the Millennium Prize problems!---might
become no harder than verifying their solutions once found.

For details about the history of P vs.\ NP,
see the superb 2022 account by Fortnow [\fort], who writes:
\begin{quote}
\noindent In those early days of P vs.\ NP, we saw NP-completeness as a
barrier---these were problems that we just couldn't solve.  As computers
and algorithms evolved, we found we could make
progress on many NP problems through a combination
of heuristics, approximation, and brute-force computing.
\end{quote}

\section{Solvability of the Navier-Stokes Equations}
For centuries, after Newton's discoveries led to the development
of mathematical physics by Euler and Lagrange and others, partial
differential equations were at the heart of science and one of the
central topics of mathematics.  The classical wave, Laplace, and heat
equations are linear, but the fundamental PDE\kern .5pt s of fluid
flow, the Navier-Stokes (NS) equations (and their predecessor the
Euler equations for flows without viscosity), are nonlinear because
of convection, which introduces the product $u\cdot \nabla u$
into the acceleration term in Newton's law.

With any PDE problem, one faces the question of whether it is
well posed.  Does a unique solution exist, and does it depend
continuously on boundary and initial data?  In much of PDE theory
and practice, these are housekeeping matters that engineers as
well as theorists deal with routinely.  The Laplace equation, for
example, is well-posed if you specify function values all around a
bounded domain, but ill-posed if you specify function values and
derivatives on a portion of the boundary.  Such facts were well
known before the formal concept of well-posedness was introduced
by Hadamard early in the 20th century.

When a PDE problem is nonlinear, understanding its behavior can be
much more challenging, and the Navier-Stokes equations have been
a source of mystery throughout their 200-year history.  The range
of behaviors of viscous fluids is enormous, with steady laminar
flow at low Reynolds numbers (nondimensionalized speeds), then
flow separation at edges, irregularity and breakdown from 2D to
3D structures at higher Reynolds numbers, and soon the appearance
of turbulence.  Flow instabilities are of ever-present importance
and are analyzed endlessly. The equations have been established
since 1846, but determining their consequences remains a perennial
project in engineering.

Given this context, it is hardly surprising that mathematicians
should be troubled that it has never been proved that the NS
equations are well-posed.  The Millennium Prize problem targets the
case of $C^\infty$ initial data in an unbounded or periodic
domain.  Despite the smoothing effect of viscosity, which becomes
all the more dominant on smaller space scales, might certain initial
conditions somehow generate singularities in finite time?

To see how natural it is to be concerned about singularity formation
in fluid flows, we need go no further than the adjacent situation of
flows that are inviscid and compressible.  Here, smooth solutions
really do develop singularities, namely {\em shock waves\/}.
This became an urgent issue in the era of jet aircraft and nuclear
explosions after the Second World War.  The hyperbolic PDE\kern
.4pt s that govern shocks do not have unique solutions.  To ensure
uniqueness, a starting point is to impose additional {\em entropy
conditions\/}, as investigated by Lax
and others in the second half of the 20th century.  It is a complex
field, and years later,
some uniqueness questions for shock waves are still not settled.

So, the problem of existence of unique solutions to the NS equations
is a natural one, even if the engineers have never been troubled
by it.  But research since 2000 has narrowed the matter.  The pursuit
of possible singularities has been intense, with heavy use of computing,
and much more is known about likely
mechanisms than before.  The work of Hou and his collaborators,
in particular, has given new shape to thoughts about how certain
initial conditions might generate a singularity [\chenhou].  But it
has proved extraordinarily hard to construct such singularities,
even while it is not known how to rule them out.  The more
we learn, the more special the possible scenarios leading to
Navier-Stokes breakdown appear to be, and the farther removed from ``wet''
fluid mechanics.  I doubt there are many these days who believe that
singularities are going to turn up in
actual flows.  The necessary initial conditions would appear to be
too contrived---and indeed, it is thought that singularity-forming
configurations, if they exist, may themselves be unstable in the
sense that small perturbations will tend to grow and shut off
the effect [\pinball,\wang].

In other words, lively research since 2000 may have rendered the problem
of well-posedness of the NS equations more theoretical than it was.
There's no predicting which way its eventual resolution may go,
but it is fascinating to speculate what may happen if it is proved
that singularities can arise.  I think that in this case, the
next scientific challenge will be not so much to modify the NS
equations to make them more physical, which might
have been the original expectation, as to
understand why those singularities have so little consequence.

\section*{Discussion}
I have argued that the Riemann Hypothesis, P vs.\ NP, and
Navier-Stokes solvability problems are all less weighty
in their direct implications now than they seemed when first
posed, even though they are fruitful and exciting as theoretical challenges.
With RH, the reason is that the
computation of trillions of Riemann zeros on the critical line has
reduced what a failure of the hypothesis would imply
quantitatively about the distribution of primes. With P vs.\ NP,
the reason is that many algorithms that are exponentially slow in
theory often prove fast in practice for reasons of typical-vs.-worst
case behavior and approximate optimality.  With Navier-Stokes,
the reason is that intense research on the problem has suggested
that flow configurations leading to singularities, if they exist,
would have to be precisely tuned to special structures and might be
in themselves unstable, diminishing their consequences for real flows.

If great challenges of mathematics tend to turn more theoretical
with the decades as they resist solution, might there be a general
explanation of this phenomenon?  One could speculate as follows.
To resolve a problem one way or another, we need to find a handle
to grab it by.  Maybe these handles have something to do
with what gives a problem, as it were, measurable consequences.
Perhaps problems that remain open for a century despite intense
efforts to solve them tend to be so smooth that they glide through
both our theory and our practice, like neutrinos, hard to catch.

I expect I will live to see the resolution of at least one of the Riemann Hypothesis,
P vs.\ NP, or the Navier-Stokes solvability problem.  This will be
thrilling, and the impact on the relevant field and indeed on
mathematics in the large will be historic.  Yet that impact will be
not so much in the direction originally imagined when the problem
was first posed, as on the rich new theories that have developed
around it subsequently.

\medskip

{\bf Acknowledgments.}
I have discussed drafts of this essay with many people, most of whom
know at least one of these problems better than I do.  Few entirely
agree with my views, but the interactions have been stimulating
and have sharpened the piece in many ways.  With gratitude, let
me mention Scott Aaronson, Folkmar Bornemann, Richard Brent,
Ernie Davis, Toby Driscoll, Ethan Epperly, Lance Fortnow, Dan
Freed, Javier Gomez-Serrano, Jonathan Goodman, Nick Gould, Tom
Hou, Russell Impagliazzo, David Jerison, Jon Keating, Don Knuth,
Kate McLoughlin, Michael Overton, Siobhan Roberts, Suchant Sachdeva, Peter Sarnak, Danny Sleator,
Endre S\"uli, and Tim Trudgian.  Two anonymous referees also improved
the presentation.

\medskip

\parindent=0pt \parskip=2pt
{\bf References}

[\zoo] S. Aaronson, Complexity Zoo, \verb|https://complex| \verb|ityzoo.net/Complexity_Zoo|.

[\apple] D. L. Applegate, R. E. Bixby, V. Chv\'atal, and W. J. Cook,
{\em The Traveling Salesman Problem: A Computational Study,} Princeton (2006). 

[\bomb]
E. Bombieri,
Problems of the millennium: The Riemann hypothesis,
{\tt https://www.claymath.org/ wp-content/uploads/2022/05/riemann.pdf}.

[\chenhou] J. Chen and T. Y. Hou, 
Singularity formation in 3D Euler equations with smooth initial data and boundary,
{\em Proc.\ Nat.\ Acad.\ Sci.,} 122 (2025).

[\cmi] Clay Mathematics Institute, The Millennium Prize Problems,
{\tt https://www.claymath.org/ millennium-problems/}.

[\cook]
S. Cook, The P versus NP problem,
{\tt https://www. claymath.org/millennium/p-vs-np/}.

[\cooktsp]
W. Cook, TSP,
{\tt https://www.math.uwaterloo .ca/tsp/optimal/} (2023).

[\ct]
J. W. Cooley and J. W. Tukey, An algorithm for the machine
calculation of complex Fourier series,
Math.\ Comput., 19 (1965), 297--301. 

[\feff]
C. L. Fefferman, Existence and smoothness of the
Navier--Stokes equation,
{\tt https://www.clay math.org/wp-content/uploads/2022/06/navier stokes.pdf\/}.

[\fort]
L. Fortnow, Fifty years of P vs.\ NP and the possibility of
the impossible, {\em Commun.\ ACM}, 65 (2022), 76--85.

[\gourdon]
X. Gourdon, Computation of zeros of the Zeta function,
{\tt http://numbers.computation.free.fr/ Constants/Miscellaneous/zetazeroscompute.
html\#Gourdon2004}.

[\mont] H. L. Montgomery and R. C. Vaughan, {\em Multiplicative
Number Theory, v.~1: Classical Theory,} Cambridge (2007).

[\pt] D. Platt and T. Trudgian, The Riemann hypothesis is true
up to $3\cdot 10^{12}$, {\em Bull.\ LMS}, 53 (2021), 792--797.

[\sarnak]
P. Sarnak, The Riemann hypothesis (video lecture, Harvard Center of
Mathematical Sciences and Applications),
{\tt https://www.youtube.com/watch? v=DtaFyE9BcXw} (2026).

[\yogiisms]
L. N. Trefethen, Inverse Yogiisms, {\em Notices AMS,} 63 (2016), 1281--1285.

[\pinball]
L. N. Trefethen, Unstable instabilities, {\em LMS Newsletter,} to appear (2026).

[\wang]
Y. Wang et al., Discovery of unstable singularities, arXiv:2509.14185v1 (2025).

[\npc]
Wikipedia, List of NP-complete problems,
\verb|https://en.wikipedia.org/wiki/List_of_NP-| \verb|complete_problems|.

[\will]
D. P. Williamson and D. B. Shmoys, {\em The Design of Approximation Algorithms\/},
Cambridge (2011).

\end{document}